\documentclass[a4paper]{amsart}
\usepackage{amsfonts,amsmath,amssymb,amscd}
\newtheorem{thm}{Theorem}[section]
\newtheorem{lem}[thm]{Lemma}
\newtheorem{prop}[thm]{Proposition}
\newtheorem{cor}[thm]{Corollary}
\newtheorem{defn}[thm]{Definition}
\newtheorem{rem}[thm]{Remark}
\def\SL{\mathop{\operatorname{SL}\nolimits}}
\def\sl{\mathop{\operatorname{\mathfrak{sl}}}\nolimits}
\def\SO{\mathop{\operatorname{SO}}\nolimits}
\def\so{\mathop{\operatorname{\mathfrak{so}}}\nolimits}
\def\O{\mathop{\operatorname{O}}\nolimits}

\def\SU{\mathop{\operatorname{SU}\nolimits}}
\def\Hom{\mathop{\operatorname{Hom}}\nolimits}
\def\ind{\mathop{\operatorname{Ind}}\nolimits}

\def\sgn{\mathop{\operatorname{sign}}\nolimits}
\def\sign{\mathop{\operatorname{sign}}\nolimits}
\def\diag{\mathop{\operatorname{diag}}\nolimits}
\def\EXP{\mathop{\operatorname{EXP}}\nolimits}
\def\vol{\mathop{\operatorname{vol}}\nolimits}
\def\tr{\mathop{\operatorname{tr}}\nolimits}
\def\Spin{\mathop{\operatorname{Spin}}\nolimits}
\def\Gal{\mathop{\operatorname{Gal}}\nolimits}
\def\Cent{\mathop{\operatorname{Cent}}\nolimits}
\def\beginpf{{\sc Proof.}}
\def\endpf{\hfill$\Box$\vskip 0.2cm}
\def\bbR{\mathbb R}
\def\bbC{\mathbb C}
\def\bbZ{\mathbb Z}
\def\bbS{\mathbb S}
\numberwithin{equation}{section}
\title[Automorphic Representations of $\SL(2,\mathbb R)$  and Quantization of Fields]{Automorphic Representations of $\SL(2,\mathbb R)$  and Quantization of Fields}
\author[Do Ngoc Diep, Do Thi Phuong Quynh]{Do Ngoc Diep${}^{1}$, Do Thi Phuong Quynh${}^{2}$}
\address{${}^{1}$ Institute of Mathematics, Vietnam Academy of Science and Technology, 18 Hoang Quoc Viet road, Cau Giay district, 10307 Hanoi, Vietnam;
{\tt Email: dndiep@math.ac.vn}}
\address{${}^{2}$ {\sc Medicine University Colllege, Thai Nguyen Universtiy, Thai Nguyen City, Vietnam;}
{\tt Email: phuongquynhtn@gmail.com}}
\begin{document}
\date{\bf Version of \today}
\maketitle
\begin{abstract} In this paper we make a clear relationship between the automorphic representations and the quantization through the Geometric Langlands Correspondence. We observe that the discrete series representation are realized in the sum of eigenspaces of Cartan generator, and then  present the automorphic representations in form of induced representations with inducing quantum bundle over a Riemann surface and then use the loop group representation construction to realize the automorphic representations. The Lanlands picture of automorphic representations is precised by using the Poisson summation formula.
\end{abstract}
\tableofcontents
\section{Introduction}
The representation theory of $\SL(2,\mathbb R)$ is well-known in the Bargman classification: every irreducible unitary/nonunitary representation is  unitarily/nonunitarily equivalent to one of representions in the five series:
\begin{enumerate}
\item the principal continuous series representations $(\pi_s,\mathcal P_s)$,
\item the discrete series representations $(\pi^\pm_k,\mathcal D_n), n\in \mathbb N, n \ne 0$,
\item the limits of discrete series representations $\mathcal (\pi^\pm_0,D_{\pm})$,
\item the complementary series representations $(\pi_s,\mathcal C_s), 0<s<1$,
\item the trivial one dimensional representation $1$,
\item the nonunitary finite dimensional representations $V_k$.
\end{enumerate}

Looking at a Fuchsian discrete subgroup $\Gamma$ of type I, i.e. 
$$ \Gamma \subseteq \SL(2,\mathbb Z), \vol(\Gamma\backslash \SL(2,\mathbb R)) < +\infty$$
one decomposes
the cuspidal parabolic part $ L^2_{cusp}(\Gamma \backslash \SL(2,\mathbb R))$ of  $L^2(\Gamma \backslash \SL(2,\mathbb R))$, that is consisting of the so called {\it automorphic representations} in subspaces of automorphic forms, of which each irreducible component with {\it multiplicity} equal to the dimension of the space of modular forms on the upper Poincar\'e half-plane 
$$\mathbb H = \{z\in \mathbb C | \Im(z) > 0\}.$$
There are a lot of studies concerning the automorphic forms and automorphic representations. Most of them realize the representations as some induced ones. Therefore some {\it clear geometric realization} of these representations should present some interests. 

In this paper, we use the {\it ideas of geometric quantization} to realize the automorphic representations in form of some Fock representations of loop algebras, see Theorem \ref{mthm}, below.
In order to do this, we first present the action of the group $\SL(2,\mathbb R)$ in the induced representation of discrete series as the action of some loop algebra/group  in the heighest/lowest weight representations.

The representations are uniquely up to equivalence defined by the character, which are  defined as some distribution function.
Beside of the main goal we make also some clear presentation of automorphic representations in this context, see Theorems \ref{thm1}, \ref{thm2}.
In the theorem \ref{thm3} we expose the trace formula using both the spectral side and geometric side.

\section{Endoscopy groups for $\SL(2,\mathbb R)$}

We introduce in this section the basic notions and notations concerning $\SL(2,\mathbb R)$, many of which are folklore or well-known but we collect all together  in order to fix a consistent system of our notations.

The unimodular group $G=\SL(2,\mathbb R)$ is the matrix group
$$\SL(2,\mathbb R) = \left\{\left.g=\begin{pmatrix} a & b\\ c & d \end{pmatrix} \right| a,b,c,d \in \mathbb R, \det g = 1  \right\}.$$
The group has finite center $\bbZ/2\bbZ$. 
It complexified group $G_\mathbb C = \SL(2,\mathbb C)$.
The unique maximal compact subgroup $K$ of $G$ is 
$$K = \left\{\left. k(\theta) = \pm\begin{pmatrix} \cos \theta & \sin\theta\\ - \sin\theta & \cos\theta\end{pmatrix} \right| \theta\in [0,2\pi)  \right\}.$$ 
The group is simple with the only nontrivial parabolic subgroup, which is minimal and therefore is a unique split Borel subgroup
$$B = \left\{\left. b= \begin{pmatrix} a & b \\ 0 & d \end{pmatrix} \right| a,b,d\in \mathbb R, ad=1 \right\}.$$ 
The Borel group $B$ is decomposed into semidirect product of its unipotent radical $N$ and a maximal split torus $T\cong \mathbb R^*_+$ and a compact subgroup $M= \{\pm 1\}$
It is well-known the Cartan decomposition $G=B\rtimes K= BK$
$$\begin{pmatrix} a & b \\ c & d\end{pmatrix} = \begin{pmatrix} y^{1/2} & x\\ 0 & y^{-1/2} \end{pmatrix}\begin{pmatrix} \pm\cos\theta & \pm\sin\theta\\ \mp\sin\theta & \pm\cos\theta \end{pmatrix}$$
and it is easy to compute
$$y = \frac{1}{c^2 + d^2}$$
$$\cos\theta = \pm y^{-1/2}d,\mbox{ or } \theta = \arccos \frac{d}{\sqrt{c^2+d^2}}$$
$$\pm y^{1/2}\sin\theta \pm y^{-1/2}x \cos\theta = b \mbox{ or  } x = \pm \frac{(b-d)}{c},$$
and the Langlands decomposition of the Borel subgroup 
$B=  M\ltimes(A \ltimes N) = MAN$.
The Lie algebra of $G=\SL(2,\mathbb R)$ is $\mathfrak g = \sl(2,\mathbb R) = \langle H,X,Y\rangle_\mathbb R$ where
$$H = \begin{pmatrix} 1 & 0\\ 0& -1\end{pmatrix},
X = \begin{pmatrix} 0 & 1\\ 0& 0\end{pmatrix},
Y = \begin{pmatrix} 0 & 0\\ 1& 0\end{pmatrix},$$
satisfying the Cartan commutation relations
$$[H,X]=2X, [H,Y]= -2Y, [X,Y] = H.$$
The Lie algebra of $A$ is $\mathfrak a = \langle H \rangle_\mathbb R$, the Lie algebra of $N$ is $\mathfrak n = \langle X \rangle_\mathbb R$.
The Lie algebra of $B$ is $$\mathfrak b = \mathfrak a \oplus \mathfrak n = \langle H,X\rangle_\mathbb R.$$
The complex Cartan subalgebra of $\mathfrak g$ is a complex subalgebra 
$$\mathfrak h = \langle H \rangle_\mathbb C \subset \mathfrak g_\mathfrak C,$$ which is concided with it normalizer in $\mathfrak g_\mathbb C$. The corresponding subgroup of $G$ such that its Lie algebra is a Cartan subalgebra, is called a Cartan subgroup.
The root system of $(\mathfrak g,\mathfrak h)$ is of rank 1 and is
$$\Sigma = \{\pm\alpha\}, \alpha = (1,-1) \in \mathbb Z(1,-1)\subset \mathbb R(1,-1).$$ There is only one positive root $\alpha = (1,-1)$, which is simple. There is a compact Cartan group $T= K = \O(2)$ The  coroot vector is
$H_\alpha = (1,-1)$ and the split Cartan subgroup of $B$ is $H = \mathbb Z/2\mathbb Z \times \mathbb R^*_+ \cong \mathbb R^*$.  

\begin{defn}
An endoscopy subgroup of $G=\SL(2,\mathbb R)$ is the connected component of identity in the centralizer of a regular semisimple element of $G$.
\end{defn}

\begin{prop} \label{endoscopy}
The only possible endoscopy subgroups of $G=\SL(2,\mathbb R)$ are itself or $\SO(2)$.
\end{prop}
\beginpf
The regular semisimple elements of $\SL(2,\mathbb R)$ are of the form $g = \diag(\lambda_1,\lambda_2)$, $\lambda_1\lambda_2 =1$. If $\lambda_1 \ne \lambda_2$, the centralizer of $g$ is the center $C(G) = \{\pm 1 \}$ of the group $\SL(2,\mathbb R)$. If $\lambda_1 = \lambda_2$ and they are real, the centralizer of $g$ is the group $\SL(2,\mathbb R)$ itself. If they are complex and their arguments are opposite, we have $g= \pm\diag(e^{i\theta},e^{-i\theta})$. In this case the connected component of identity of the centralizer is $\SO(2)$. The endoscopy groups  $\SL(2,\mathbb R)$ itself or the center $\{\pm 1\}$ are trivial and the only nontrivial endoscopy group is $\SO(2)$.
\endpf

\section{Automorphic Representations}
In this section we make clear the construction of automorphic representations.

The following lemma is well-known.
\begin{lem}
There is a one-to-one correspondence between any irreducible $2n-1$ dimensional representations of $\SO(3)$ and the $n$ dimensional representations of $\SO(2)$.
\end{lem}
\beginpf
There is a short exact sequence
\begin{equation} \CD
1 @>>> \{\pm I\} @>>> \SU(2) @>>> \SO(3) @>>> 1
\endCD\end{equation}
The characters of the $2n-1, n=1,2,\dots$ dimensional representation of $\SO(3)$ is 
\begin{equation}
\chi_{2n-1}(k(\theta)) = \frac{\sin (2n-1)\theta}{\sin \theta}
\end{equation}
where $$k(\theta) = \begin{pmatrix}
\cos\theta & \sin\theta\\  -\sin\theta & \cos\theta \end{pmatrix}$$
\endpf

It is well-known the Iwasawa decomposition $ANK$: every element $g= \begin{pmatrix} a & b\\ c & d \end{pmatrix}$ has a unique decomposition of form 
\begin{equation}\label{iwasawa}
\begin{pmatrix} a & b\\ c & d \end{pmatrix} =
\begin{pmatrix} y^{1/2} & 0\\  0 & y^{-1/2} \end{pmatrix}
\begin{pmatrix} 1 & y^{-1/2}x\\ 0 & 1 \end{pmatrix}
\begin{pmatrix} \cos\theta & -\sin\theta\\ \sin\theta & \cos\theta \end{pmatrix}
\end{equation}
of $\SL_2(\mathbb R)$, where
\begin{equation}
N = \left\{ \begin{pmatrix} 1 & x\\ 0 & 1 \end{pmatrix} \right\} \end{equation}
 is the unipotent radical of the Borel subgroup,
\begin{equation}
A = \left\{ \begin{pmatrix} y^{1/2} & 0\\  0 & y^{-1/2} \end{pmatrix}\right\}
\end{equation}
is the split torus in the Cartan subgroup, and
\begin{equation}
K = \left\{\begin{pmatrix} \cos\theta & \sin\theta\\ -\sin\theta & \cos\theta \end{pmatrix}\right\} \end{equation}
is the maximal compact subgroup.
\vskip .2cm

To each modular form $f\in\mathcal S_k(\Gamma)$ of weight $k$ on the Poincar\'e plane $\mathbb H = \SL(2,\mathbb R)/\SO(2)$, we associate the automorphic form $\varphi_f \in \mathcal A_k(\SL(2,\mathbb R))$
\begin{equation}
\varphi_f(\begin{pmatrix} a & b\\ c & d \end{pmatrix} ) = y^{k/2}e^{ik\theta} f(x+iy),
\end{equation}
where $x,y,\theta$ are from the Iwasawa decomposition (\ref{iwasawa}).

The discrete series representations are realized on the function on $L^2(\mathbb H,\mu_k)$, $\mu_k = y^k\frac{dxdy}{y^2}$ of weight $k$ modular form by the formula
\begin{equation}
\pi_k(\begin{pmatrix} a & b\\ c & d \end{pmatrix})f(z) = (cz+d)^{-k}f(z). \end{equation}
Denote by $\mathcal D_k$ the discrete series representation of weight $k$.
The cuspidal automorphic representations are realized in the space $L^2_{cusp}(\Gamma\backslash \SL(2,\mathbb R))$ of automorphic forms

\begin{thm}\label{thm1}
The set of interwining $\SL(2, \mathbb R)$ homomorphisms  from the set of discrete series representations to the set $L^2_{cusp}(\Gamma \backslash\SL(2, \mathbb R))$ of the  automorphic representations of 
$\SL(2, \mathbb R)$ is equal to the set $S_k(\Gamma)$ of modular forms, i.e.
\begin{equation} 
\Hom_{\SL(2, \mathbb R)}(\mathcal D_k, L^2_{cusp}(\Gamma \backslash \SL(2, \mathbb R)) = \mathcal S_k(\Gamma)
\end{equation}
\end{thm}
\beginpf
The theorem is well-known in literature, see for example A. Borel \cite{borel}. Starting from some intertwining operator $A\in \Hom_{\SL(2, \mathbb R)}(\mathcal D_k,L^2_{cusp}(\Gamma \backslash\SL(2, \mathbb R)) $ we may construct the $L$-series which is an element in $S_k(\Gamma)$; and conversely, starting from some a modular form $f\in S_k(\Gamma)$ we construct the corresponding $L$-series $L_f$. There exists a unique interwining operator $A$ such that the $L$-series of which is equal to $L_f$. 
\endpf

\subsection{Geometric Langlands Correspondence}
The general Geometric Langlands Conjecture was stated by V. Drinfel'd and was then proven by E. Frenkel, D. Gaitsgory and Vilonen \cite{vilonenetc} and became the Geometric Langlands Correspondence (GLC). 
We will start to specify the general GLC in our particular case of group $\SL_2(\mathbb R)$.

\begin{thm}[Geometric Langlands Correspondence]\label{thm2}
There is a bijection 
\begin{equation}
[\pi_1(\Sigma), \SO(3)] \longleftrightarrow \mathcal A(\SL(2,\mathbb R)) 
\end{equation}
between the set  of equivalence classes of representation of the fundamental group $\pi_1(\Sigma)$ of the Riemannian surface $\Sigma = \Gamma\backslash\SL(2,\bbR)/\SO(2)$ in $\SO(3)$ and the set $\mathcal A(\SL(2, \mathbb R))$ of equivalence classes of automorphic representations of $\SL(2, \mathbb R)$.
\end{thm}
\beginpf The theorem was proven in the general context of a reductive group in the works of \cite{vilonenetc}. For the particular case of $\SL(2,\mathbb R)$ many things are simplified, what we want to point out here.

The main idea to prove the theorem is consisting of the following ingredients:
\begin{enumerate}
\item
It is well-known that every element of $\SO(3)$ is conjugate with some element of the maximal torus subgroup $\SO(2)$. Therefore the set of homomorphism from $\pi_1(\Sigma)$ into $\SO(3)$ is the same as the set of character of the Borel subgroup (minimal parabolic subgroup) from which the discrete series representations are induced. 
\begin{lem}
There is one-to-one correspondencee between the conguacy classes of $\SO(3)$ in itself and the inducing character for the discrete series representations.
\end{lem}
\item
Moreover, the $\Gamma$ invariance condition is the same as the
 condition to be extended from the local character of some automorphic component. The following two lemmas are more or less known \cite{bailly}.
\begin{lem}\label{conjcl}
Every representation of $\pi_1(\Sigma)$ in $\SO(3)$ is defined by a system of conjugacy classes in $\SO(3)$, one per generator of $\pi_1(\Sigma)$ 
\end{lem}

\begin{lem}
Every system of conjugacy classes in the previous lemma \ref{conjcl} defines a unique modular form on $\mathbb H$ and hence a unique automorphic representation of $\SL_2(\mathbb R)$.
\end{lem}

\end{enumerate}
\endpf

\subsection{Geometric Quantization}
The idea of realizing the automorphic representations of reductive Lie groups was done in \cite{diep1}. In this section we show the concrete computation for the case of $\SL(2, \mathbb R)$.

\begin{thm}\label{mthm}
The automorphic representations are obtained from the quantization procedure of fields based on geometric Langlands correpondence.
\end{thm}
\beginpf
The discrete series representations can be realized through the geometric quantization as follows. 
\begin{enumerate}
\item
The representation space of discrete series representation is consisting of square-intergrable holomorphic functions. 
\begin{equation} \label{convseries}
f(z) = \sum_{n=0}^\infty c_nz^n; \quad \sum_{n=0}^\infty |c_n|^2 < \infty \end{equation} 
\begin{lem}
The Hardy space of holomorphic quare-integrable functions can be realized as the expnential Fock space of the standard representation of $\SO(2)$ in $\mathbb C$.
\end{lem}
Indeed every module-square convergence series of type ref{convseries} can be express as some element
\begin{equation} \label{convseries}
f(z) = \sum_{n=0}^\infty n!c_n\frac{z^n}{n!}; \quad \sum_{n=0}^\infty |c_n|^2 < \infty \end{equation} 
 in the exponential vector space 
\begin{equation}\label{exp}
\EXP \mathbb C = \bigoplus_{n=0}^\infty \frac{C^{\otimes n}}{n!},\quad \mathbb C^{\otimes n} \cong \mathbb C
\end{equation}
\item
The Lie algebra $\sl_2(\mathbb R)$ with 3 generators $\sl_2(\mathbb R) = \langle H, X, Y \rangle_\mathbb R$ in the induced representations of discrete series act through the action of one-parameters subgroups
\begin{equation}
g_3(t) = \exp (tH) = \begin{pmatrix} e^t & 0\\ 0 & e^{-t} \end{pmatrix}
\end{equation}
$$g_1(t) = \exp(t(X-Y)) = \begin{pmatrix} \cos t & \sin t\\ -\sin t & \cos t \end{pmatrix}$$
$$g_2(t) = \exp(t(X+Y)) = \begin{pmatrix} \cosh t & \sinh t\\ \sinh t & \cosh t \end{pmatrix}$$
Under the representation $\pi^{\pm}_n$ we have
\begin{equation} 
\pi^{\pm}_n(g_k(t)) = e^{it\hat U_k}
\end{equation}
and we have the relations
\begin{equation}
i\hat U_3 = \hat H = 2z\frac{\partial}{\partial z} +(n+1),
\end{equation}
$$i\hat U_1 = \hat X - \hat Y = -(1+z^2)\frac{\partial}{\partial z} - (n+1)z,$$
$$i\hat U_2 = \hat X + \hat Y = (1-z^2)\frac{\partial}{\partial z} - (n+1)z .$$

In this representation, the action of the element 
\begin{equation}
\Delta = \frac{-1}{4}(\hat U_1^2 - \hat U_2^2 - \hat U_3^2 )=
\end{equation}

 $$=\frac{-1}{4}\left((\hat X - \hat Y)^2 - (\hat X + \hat Y)^2 - (\hat H)^2\right) = \frac{1}{4} (\hat H^2  + 4\hat X\hat Y)=$$
$$= (z -\bar z)^2\frac{\partial^2}{\partial z\partial \bar{z}}= -y^2\left( \frac{\partial^2}{\partial x^2} + \frac{\partial^2}{\partial y} \right)$$

This action can be represented as the action of the loop algebra over the Riemann surface $\Sigma$ with values in $\SO(2)$, the elements of the loop algebra are presented in form of some formal/conformal Laurent series with values in the corresponding Lie algebra of form of connection appeared in the construction of induced representations, i.e.
\begin{equation} T(z) = \sum_{n=-\infty}^\infty c_nz^n, c_n\in \so(2), z\in\Sigma . \end{equation}
In our case the Lie algebra $\so(2)$ is one dimensional and we have all number coefficients $c_n, n\in \mathbb Z$.
\begin{lem}
The decompositon \ref{exp}  of $\EXP \mathbb C$ presents the weight decomposition of $\sl_2(\mathbb R)$ in which $H$ keeps each component, $X$ acting as creating operator and $Y$ is acting as some annihilating operator. 
\end{lem}
\item The lowest weight representations of the loop algebras are realized through the lowest weight representations of the Virassoro algebra as follows. Let us consider the generators
\begin{equation} L_n = \int_\Sigma z^{n+1} T(z)dz \end{equation}
for any element $$T(.): \Sigma \to \SO(2)$$ from the loop algebra presentation are realized in the Fock space of the standard representation.
These generators satisfy the Virasoro algebra relations
\begin{equation}
[L_m,L_n] = (n-m)L_{m+n} + \delta_{n,-m}\frac{n(n^2-1)}{12}L_0,
\end{equation}
where in an irreducible representation, $Z=cL_0= cI$ is the central charge element.
\end{enumerate}
The proof of Theorem \ref{thm2} is therefore achieved.
\endpf

\section{Arthur-Selberg Trace Formula}
\subsection{Trace Formula}
Let us review the trace formula due to Selberg and J. Arthur.
Remind that by $\mathbb H = \SL(2,\mathbb R) \backslash \SO(2)$ denote the upper Poincar\'e half-plane, $$\mathbb H = \{ z= x+iy \in \mathbb C | \Im z = y > 0 \},$$ 
$\Gamma$ denote a Langlands type discrete subgroup of finite type with finite number of cusps $\kappa_1,\dots,\kappa_h$. Let $$\Gamma_0 =\left\{\left. \gamma = \begin{pmatrix}
1 & b\\ 0 & 1 \end{pmatrix} \right| b \in \mathbb Z  \right\} \subset \SL(2,\mathbb Z)$$ and $$\Gamma_i= \{\sigma \in \Gamma | \sigma \kappa_i = \kappa_i  \} \subset \SL(2,\mathbb Z)$$ and $\sigma_i\in \SL(2,\mathbb R), i=\overline{1,h}$ are such that  $\sigma_i\Gamma_i\sigma_i^{-1} = \Gamma_0$.
$\sigma_i\Gamma_i\sigma_i = \Gamma_0$,
$\mathfrak H = L^2(\Gamma\backslash G)$ denote the space of quare-integrable functions on $G$ on which there is a natural rugular representation $R$ of $G$ by formula
$$[R(\begin{pmatrix} a & b \\ c & d \end{pmatrix})f](z) = f(\frac{az+b}{cz+d}), z \in \mathbb H.$$ In particular, the unipotent radical $N= \left\{\left. \begin{pmatrix} 1 & x\\ 0 & 1 \end{pmatrix} \right| x \in \mathbb R  \right\}$ of $B$ is acting by the translation on varable $z$ to $z +x$. For any funciton $\psi: N \backslash \mathbb H \to \mathbb C$ of variable $y$ decreasing fast enough when $y$ approaches to 0 or $\infty$,  one defines the incomplete $\theta$-series 
$$\theta_{t,\psi}(z) = \sum_{\sigma\in \Gamma_i\backslash \Gamma}  \psi(\sigma_i^{-1}\sigma z)$$
which is certainly of class $L^2(N\cap \Gamma\backslash \mathbb H)$. Denote by $$\Theta = \langle \theta_{t,\psi} | \forall \psi, t \rangle\subset L^2(\Gamma \backslash H)$$ the space of incomplete $\theta$-series

It is well-known that the orthogonal complement $\Theta^\perp$ of $\Theta$ in $L^2(\Gamma \backslash \mathbb H)$ is isomorphic to the space $\mathfrak H_0$ of cuspidal parabolic forms or in other words, of automorphic forms with zero Fourier constant  terms of  automorphic representations, $\mathfrak H = \Theta \oplus \mathfrak H_0$.

Consider the Hecke algebra $\mathcal H(\SL(2,\mathbb R))$ of all convolution Hecke operators of form as follows. Let $F : \mathbb H \to \mathbb C$ be a function $K$-invariant with respect to transformations of form $z \mapsto \gamma z$, for all $\gamma \in K$. Such a function is uniquely defined by a function on $K \backslash G \slash K$, or a so called spherical function $F$ on $G$ which is left and right invariant by $K$. Under convolution these functions provide the Hecke algebra of Hecke operators by convolution with functions in representations. All the Hecke operators have kernel as follows.
For any automorphic function $f$. $f(\gamma z) = f(z), \forall \gamma \in \Gamma$,
$$(F*f)(z) = \int_\mathbb H F({g'}^{-1}g)f(g')dg' = \int_\mathbb H k(z,z')f(z')dz'= \int_{\Gamma\backslash \mathbb H} \sum_{\gamma\in \Gamma} k(z,\gamma z')f(z')dz',$$
$$= \int_{\Gamma\backslash \mathbb H} K(z,z')f(z')dz',$$
 where $z = gi, z' g'i, i = \sqrt{-1}$, $$K(z,z')= \sum_{\gamma\in \Gamma}k(z,\gamma z').$$ Denote the sum of kernel over cups  by  
$$H(z,z') = \sum_{i=1}^h H_i(z,z') = \sum_{i=1}^h\sum_{\gamma_i\in \Gamma_i\backslash \Gamma} \int_{-\infty}^{+\infty} K(z,\gamma_i n(x)\gamma_i^{-1}\gamma z')dx$$
The Hecke operator with kernel $K(z,z')$ has the same spectrum as the operator with kernel $K^*(z,z') = K(z,z') - H(z,z')$. The kernel $K^*(z,z')$ are bounded and the fundamental domain $\Gamma\backslash \mathbb H$ has finite volume, therefore the kernels $K^*(z,z')$ are of class $L^2$ on $\mathcal D \times \mathcal D, \mathcal D = \Gamma\backslash \mathbb H$, and all the Hecke operators are compact operators. The Hecke operators keeps each irreducible components of $\Theta^\perp$ invariant and therefore are scalar on each automorphic representation.
On each irreducible component, the Laplace operators has also a fixed eigenvalue 
$$\Delta f = \lambda f, \lambda = \frac{s(s-1)}{4}, \Delta = -y^2\left(\frac{\partial^2}{\partial x^2} + \frac{\partial^2}{\partial y^2}  \right)$$

One deduces therefore the theorem of spectral decomposition for the discrete part of the regular representation.

\begin{thm}[Spectral decomposition] In the induced representation space of $\ind_B^G\chi_{\lambda,\varepsilon}$, choose $$H_n = \left\{f \in H | \pi(\begin{pmatrix} \cos\theta & \sin\theta\\ -\sin\theta & \cos\theta \end{pmatrix})f = e^{in\theta}f \right\}.$$. They are all of dimension 1 and $$H = \bigoplus_{n\in \mathbb Z} H_n.$$
The discrete part $R|_{L^2_{cusp}(\Gamma\backslash \SL(2,\mathbb R)}$ of the regular representation can be decomposed as the sum of the discrete series representations $\pi_n^\pm$ in spaces $$D_{s+1}^+=\bigoplus_{n\equiv\varepsilon \mod 2, n \geq m} H_n$$ or $$D_{s+1}^-=\bigoplus_{n\equiv\varepsilon \mod 2, n \leq -m}H_n,$$ $s\in \mathbb Z, s>0$ and $s+1 \equiv \varepsilon \mod 2$ and there exists $m\in \mathbb Z, m=s+1, m>0$, induced from $\chi_{i\lambda,\varepsilon} = |a|^{i\lambda}(\sgn a)^\varepsilon$ and limits of principal series representations $\pi_0^\pm$ in $\mathcal D_1^+$ or $\mathcal D_1^-$ as two component of the representation $\pi^{0,1} = \ind_B^G\chi_{0,1}$,   induced from the character $\chi_{0,1}$. 
\end{thm}
\begin{rem}
In the spaces $\bigoplus_{-m<n<m} H_n$ of dimension $2m-1$ the finite dimensional representations $V_m$ are realized. 
\end{rem}

\begin{cor}
For any function $\varphi$ of class $C^\infty_c(G)$, the operator $\pi^\pm_n(\varphi)$ is of trace class and is a distribution denoted by $\Theta^\pm_n$ (following Harish-Chandra) which are uniquely defined by their restriction to the maximal compact subgroup $K = \SO(2)$ and
$$\tr R(\varphi) = \sum_{n \in \mathbb Z, n\geq 0, \pm} m(\pi_n^\pm)\Theta^\pm_n(\varphi) ,$$ with multiplicities $m(\pi^\pm_n)$.
\end{cor}

\subsection{Stable Trace Formula}

The Galois group $\Gal(\mathbb C/\mathbb R) = \mathbb Z_2$ of the complex field $\mathbb C$ is acting on the discrete series representation by character $\kappa(\sigma) = \pm 1$. Therefore the sum of characters can be rewrite as some sum over stable classes of characters.
$$\tr R(f) = \sum_{n=1}^\infty \sum_{\varepsilon = \pm 1} (\Theta_n^+(f) - \Theta_n^-(f)).$$
\begin{rem}The stable trace formula is uniquely defined by its restriction to the maximal compact subgroup $K=\SO(2)$ and is
$$S\Theta_n = \frac{\Theta_{n}^+ - \Theta_n^-}{e^{i\theta} - e^{-i\theta}} = -2i\sin\theta(e^{in\theta}-e^{-in\theta}).$$
\end{rem}

\section{Endoscopy}
\subsection{Stable orbital integral}
Let us remind that the \textit{orbital integral} is defined as
$$\mathcal{O}_\gamma(f) = \int_{G_\gamma\backslash G} f(x^{-1}\gamma x)d\dot x$$ The main idea of computation of endoscopy transfer was explained in \cite{labesse}, we make it here in more detail to clarify some points.  

The complex Weyl group is isomorphic to $\mathfrak S_2$ while the real Weyl group
is isomorphic to $\mathfrak S_1$ . The set of conjugacy classes inside a strongly regular
stable elliptic conjugacy class is in bijection with the pointed set
$\mathfrak S_2 /\mathfrak S_1= \mathfrak S_2$
 that can be viewed as a sub-pointed-set of the group
$\mathfrak E(\mathbb R, T, G) = (Z_2 )$
We shall denote by $\mathfrak K(\mathbb R, T, G)\cong \bbZ_2$ its Pontryagin dual.

Consider $\kappa \ne 1$ in $\mathfrak K(\mathbb R, T, G)$ such that $\kappa(H) = -1$. 
 Such a $\kappa$ is unique.
The endoscopic group $H$ one associates to $\kappa$ is isomorphic to
$\SO(2)$
the positive root of $\mathfrak h$ in $H$ (for a compatible order) being $\alpha = \rho$

Let $f_\mu$ be a pseudo-coefficient for the discrete series representation $\pi_\mu$ then the \textit{$\kappa$-orbital integral} of a  regular element $\gamma$
in $T (R)$ is given by
$$\mathcal O^\kappa_\gamma(f_\mu) =\int_{G_\gamma\backslash G} \kappa(x) f_\mu(x^{-1}\gamma x) d\dot x = \sum_{\sign(w) =1} \kappa(w) \Theta^G_\mu(\gamma^{-1}_w)  = \sum_{\sign(w) =1} \kappa(w)\Theta_{w\mu}(\gamma^{-1}),$$ because there is a natural bijection between the left coset classes and the right coset classes.

\subsection{Endoscopic transfer}
We make details for the guides of J.-P. Labesse \cite{labesse}.

\textsl{The simplest case} is the case when $\gamma = \diag(a, a^{-1})$. In this case, because of Iwasawa decomposition $x=auk$, and the $K$-bivariance,  the orbital integral is
$$\mathcal O_\gamma(f) = \int_{G_\gamma\backslash G} f(x^{-1}\gamma x)dx = \int_U f(u^{-1}\gamma u)du$$ $$=\int_\mathbb R f(\begin{pmatrix} 1 & x\\ 0 & 1\end{pmatrix}^{-1}\begin{pmatrix} a & 0 \\ 0 & a^{-1} \end{pmatrix} \begin{pmatrix} 1 & x\\ 0 & 1\end{pmatrix})dx$$ 
$$=\int_\mathbb R f(\begin{pmatrix} 1 & -x\\ 0 & 1\end{pmatrix}\begin{pmatrix} a & 0 \\ 0 & a^{-1} \end{pmatrix} \begin{pmatrix} 1 & x\\ 0 & 1\end{pmatrix})dx$$ $$= \int_\bbR f(\begin{pmatrix} a & (a-a^{-1})x\\ 0 & a^{-1}\end{pmatrix})dx = |a-a^{-1}|^{-1} \mathcal O_\gamma(f).$$
The integral is abosolutely and uniformly convergent and therefore is smooth function of $a\in \bbR^*_+$. Therefore the function $$f^H(\gamma) = \Delta(\gamma)\mathcal O_\gamma(f), \quad \Delta(\gamma) = |a- a^{-1}|$$ is a smooth function on the endoscopic group $H= \bbR^*$. 

\textsl{The second case} is the case where $\gamma = k_\theta = \begin{pmatrix} \cos\theta & \sin\theta\\ -\sin\theta & \cos\theta \end{pmatrix}$. We have again, $x= auk$ and
$$\mathcal O_{k(\theta)}(f) = \int_{G_{k(\theta)}\backslash G} f(k^{-1}u^{-1}a^{-1}k(\theta)auk)dx\frac{dy}{|y|}d\theta$$
$$ = \int_{G_{k(\theta)}\backslash G} f(u^{-1}a^{-1}k(\theta)au)dx\frac{dy}{|y|}d\theta$$
$$ = \int_{G_{k(\theta)}\backslash G} f(\begin{pmatrix} 1 & -x \\ 0 & 1 \end{pmatrix}\begin{pmatrix} a^{-1} & 0\\ 0 & a \end{pmatrix} \begin{pmatrix} \cos\theta & \sin\theta\\ -\sin\theta & \cos\theta \end{pmatrix}\begin{pmatrix} a^{-1} & 0\\ 0 & a \end{pmatrix}\begin{pmatrix} 1 & x \\ 0 & 1 \end{pmatrix})dx\frac{dy}{|y|}d\theta$$
$$= \int_1^\infty f(\begin{pmatrix}\cos\theta & t\sin\theta\\ -t^{-1}\sin\theta & \cos\theta \end{pmatrix})|t-t^{-1}|\frac{dt}{t}$$ $$= \int_0^{+\infty}\sgn(t-1) f(\begin{pmatrix}\cos\theta & t\sin\theta\\ -t^{-1}\sin\theta & \cos\theta \end{pmatrix})dt.$$ When $f$ is an element of the Hecke algebra, i.e. $f$ is of class $C^\infty_0(G)$ and is $K$-bi-invariant, the integral is converging absolutely and uniformly. Therefore the result is a function $F(\sin\theta)$. 
The function $f$ has compact support, then the integral is well convergent at $+\infty$. At the another point $0$, we develope the function $F$ into the Tayor-Lagrange of the first order with respect to $\lambda = \sin\theta \to 0$
$$F(\lambda) = A(\lambda) + \lambda B(\lambda),$$ where $A(\lambda) = F(0)$ and $B(\lambda)$ is the error-correction term $F'(\tau)$ at some intermediate value $\tau, 0 \leq \tau \leq t$.
Remark that 
$$ \begin{pmatrix} \sqrt{1-\lambda^2} & t\lambda\\ -t^{-1}\lambda & \sqrt{1-\lambda^2} \end{pmatrix}=\begin{pmatrix} t^{1/2}& 0\\ 0 &t^{-1/2} \end{pmatrix} \begin{pmatrix} \sqrt{1-\lambda^2} & \lambda\\ \lambda & \sqrt{1-\lambda^2} \end{pmatrix}\begin{pmatrix} t^{-1/2} & 0\\ 0 &t^{1/2} \end{pmatrix} $$ we have
$$B= \frac{dF(\tau)}{d\lambda} =\frac{d}{d\lambda} \left. \int_0^{+\infty} \sgn(t-1)f(\begin{pmatrix} \sqrt{1-\lambda^2} & t\lambda\\ -t^{-1}\lambda & \sqrt{1-\lambda^2} \end{pmatrix})dt \right|_{t=\tau}$$
$$=  \int_0^{+\infty} \sgn(t-1)g(\begin{pmatrix} \sqrt{1-\lambda^2} & t\lambda\\ -t^{-1}\lambda & \sqrt{1-\lambda^2} \end{pmatrix})\frac{dt}{t}, $$
where $g\in C^\infty_c(N)$
and $g(\lambda) \cong O(-t^{-1}\lambda)^{-1}.$ $B$ is of logarithmic growth and $$B(\lambda) \cong \ln(|\lambda|^{-1})g(1)$$ up to constant term,  and therefore is contimuous. 
$$A = F(0) = |\lambda|^{-1} \int_0^\infty f(\begin{pmatrix}1 & \sgn(\lambda)u\\ 0 & 1 \end{pmatrix}) du - 2f(I_2) + o(\lambda)$$
Hence the functions 
$$G(\lambda) = |\lambda|(F(\lambda) + F(\lambda)),$$
$$H(\lambda) = \lambda(F(\lambda) - F(-\lambda))$$ have the Fourier decomposition
$$G(\lambda) = \sum_{n=0}^N (a_n|\lambda|^{-1} + b_n)\lambda^{2n} + o(\lambda^{2N})$$
$$H(\lambda) = \sum_{n=0}^N h_n\lambda^{2n} + o(\lambda^{2N})$$
Summarizing the discussion, we have that in the case of $\gamma = k(\theta)$, there exists also a continuous function $f^H$ such that $$f^H(\gamma) = \Delta(\gamma) (\mathcal O_\gamma(f) - \mathcal O_{w\gamma}(f))=\Delta(k(\theta))\mathcal {SO}_\gamma(f), $$ where $\Delta(k(\theta)) = -2i\sin\theta$.

\begin{thm}
There is a natural function $\varepsilon : \Pi \to \pm 1$ such that in the Grothendieck group of discrete series representation ring, $$\sigma_G = \sum_{\pi\in \Pi} \varepsilon(\pi)\pi,$$ the map $\sigma \mapsto \sigma_G$ is dual to the map of geometric transfer, that for any $f$  on $G$, there is a unique $f^H$ on $H$
$$\tr \sigma_G(f) = \tr\sigma(f^H).$$
\end{thm}
\textsc{Proof.}
There is a natural bijection $\Pi_\mu \cong \mathfrak D(\mathbb R,H,G)$, we get a pairing
$$\langle .,.\rangle : \Pi_\mu \times \mathfrak k(\mathbb R,H,G)\to \mathbb C.$$ Therefore we have
$$\tr \Sigma_\nu(f^H) = \sum_{\pi\in \Pi_\Sigma} \langle s,\pi\rangle \tr \pi(f).$$
\hfill$\Box$

Suppose given a complete set of endoscopic groups $H = \mathbb S^1 \times \mathbb S^1 \times \{\pm 1\}$ or $\SL(2,\mathbb R) \times \{\pm 1\}$.  For each group, there is a natural inclusion
$$\eta: {}^LH \hookrightarrow {}^LG$$

Let $\varphi: DW_\bbR  \to {}^LG$ be the Langlands parameter, i.e. a homomorphism from the Weil-Deligne group $DW_\bbR = W_\bbR\ltimes \bbR^*_+$ the Langlands dual group,
$\bbS_\varphi$ be  the set of conjugacy classes of Langlands parameters modulo the connected component of identity map. For any $s\in 
\bbS_\varphi$,
$\check{H}_s = \Cent(s,\check G)^\circ$ the connected component of the centralizer of $s\in\bbS_\varphi$ we have $\check{H}_s$ is conjugate with $H$. 
Following the D. Shelstad pairing $$\langle s, \pi\rangle : \bbS_\varphi \times \Pi(\varphi) \to \bbC$$
$$\varepsilon(\pi) = c(s)\langle s,\pi\rangle.$$
Therefore, the relation
$$\sum_{\sigma\in \Sigma_s} \tr \sigma(f^H) = \sum_{\pi\in\Pi} \varepsilon(\pi) \tr \pi(f)$$
can be rewritten as
$$\widetilde{\Sigma}_s(f^H) = \sum_{s\in\Pi} \langle s,\pi\rangle \tr\pi(f)$$
and
$$\widetilde{\Sigma}_s(f^H) = c(s)^{-1}\sum_{\sigma\in\widetilde{\Sigma}_s} \tr \sigma(f^H) .$$
We arrive, finally to the result
\begin{thm}
$$\tr \pi(f) = \frac{1}{\#\bbS_\varphi}\sum_{s\in\bbS_\varphi} \langle s,\pi\rangle\widetilde{\Sigma}_s(\check{f}^H).$$
\end{thm}

\section{Poisson Summation Formula}
In the Langlands picture of the trace formula, the trace of the restriction of the regular representation on the cuspidal parabolic part is the coincidence of of the spectral side and the geometric side.
\begin{equation}
\sum_{\pi} m(\pi) \hat{f}(\pi) = \sum_{\gamma\in \Gamma\cap H} a^G_\gamma \hat{f}(\gamma) \end{equation}
Let us do this in  more details.

\subsection{Geometric side of the trace formula}
\begin{thm} \label{thm61}
The trace formula for the regular representation of $\SL_2(\mathbb R)$ in the space of cusp forms is deccomposed into the sum of traces of automorphic representations with finite multiplicities is transfered into the modified Poisson summation formula
\begin{equation}
\sum_{\gamma\in \Gamma\cap H} \varepsilon(\gamma)\mathcal O_\gamma(f) = \sum_{\gamma\in\Gamma\cap H} \varepsilon(\gamma)\vol(\Gamma\cap H) \int_{H\backslash G} f(x^{-1}\gamma x)dx 
\end{equation}

\end{thm}
\beginpf
It is easy to see that the restriction of the Laplace operator $\Delta$ on the Cartan subgroup $H$ is elliptic and therefore the Cauchy problem for the other variables has a unique solution. The solution is the trace formula for the cuspidal parabolic part of the regular representation.
\begin{equation} \tr R(f)|_{ L^2_{cusp}(\Gamma\backslash G)} = \sum_{\gamma\in \Gamma\cap H} \varepsilon(\gamma)\mathcal O_\gamma(f)  \end{equation}
From another side  we have
\begin{equation}\tr R(f)|_{ L^2_{cusp}(\Gamma\backslash G)}=
\sum_{\gamma\in\Gamma\cap H}  \varepsilon(\gamma)\vol(\Gamma\cap H) \int_{H\backslash G} f(x^{-1}\gamma x)dx$$ $$=\sum_{\gamma\in \Gamma\cap H} \varepsilon(\gamma)\vol(\Gamma\cap H)\mathcal {SO}_\gamma(f^H).
\end{equation}
\endpf

\subsection{Spectral side of the trace formula}
The following result is well-known, see e.g. \cite{ggps}
\begin{thm}[Gelfand - Graev - Piateski-Shapiro] \label{GGPS}
For any compactly supported function $f\in C^\infty_0(\SL_2(\mathbb R))$ the operator $R(f)|_{L^2_{cusp}(\Gamma\backslash \SL_2(\mathbb R))}$ is of trace class and each irreducible component is of finite multiplicity
\begin{equation}
R(f)|_{L^2_{cusp}(\Gamma\backslash \SL_2(\mathbb R))} = \sum_{\pi\in \mathcal A(\SL_2(\mathbb R))} m(\pi) \pi(f), 
\end{equation}
where $m(\pi) = \dim_{\mathbb C}\Hom_{ \SL_2(\mathbb R)}(\mathcal D_k, L^2_{cusp}(\Gamma\backslash \SL_2(\mathbb R)))$.
\end{thm}

\subsection{Poisson Summation formula}
Therefore, following the Poisson summation formula we have the equlity of the both sides.  Let us denote by $\chi_k$ the character of $\SO(2)$ that induce  the discrete series representation $\mathcal D_k$ of  $\SL_2(\mathbb R)$. In our case the Cartan subgroup is $\SO(2)$ and we have the ordinary 
\begin{lem} \label{PS}
De note the universal covering of $H= \SO(2)$ by $\tilde H = \Spin(2)$
\begin{equation} 
\sum_{\gamma\in \Gamma\cap \tilde H}\delta(x+\gamma) = \sum_{(\varepsilon,k)\in\mathbb Z_2 \times \mathbb Z=\hat H}  (\sign x)^\varepsilon\chi_k(x)
\end{equation}
\end{lem}
\beginpf It is the same as ordinary Poisson summation, lifted to the universal covering $\tilde H$ of $\SO(2)$. 
\endpf

\begin{thm} \label{thm3}
The trace $\tr R(f)$ of the restriction of the regular representation on the cuspidal parabolic part $ L^2_{cusp}(G)$ is computed by the formulas
\begin{equation} \sum_{\pi\in \mathcal{A}(G)} m(\pi) \pi(f) = \sum_{\gamma\in \Gamma\cap \tilde H\atop \varepsilon(\gamma) = \pm 1}\varepsilon(\gamma) \vol(\Gamma\cap \tilde H)\mathcal O_\gamma(f)=\sum_{\gamma\in \Gamma\cap H} \vol(\Gamma\cap H)\mathcal {SO}_\gamma(f^H)\end{equation}
\end{thm}
\beginpf
The proof a combination of the previous theorems \ref{thm61}, \ref{GGPS} and \ref{PS} therefore is complete.
\endpf


\end{document}